\documentclass[reqno]{amsart}
\usepackage[english]{babel}
\usepackage{amssymb}
\usepackage{verbatim,graphicx}
\usepackage[T1,T5]{fontenc}

\theoremstyle{definition}

\theoremstyle{remark}

\numberwithin{equation}{section}

\newcommand{\dom}{\mathrm{dom}}

\newcommand{\pqbin}[2]{\genfrac{(}{)}{0pt}{0}{#1}{#2}_{\!\!p,q}}

\newcommand{\qbin}[3]{\genfrac{(}{)}{0pt}{0}{#1}{#2}_{\!\!#3}}

\begin{document}

\author{P. {\AA}hag} \address{Department
  of Mathematics and Mathematical Statistics, Ume{\aa} University,
  SE-901 87 Ume{\aa}, Sweden}
 \email{per.ahag@umu.se}

\author{R. Czy{\.z}}\address{Faculty of Mathematics and Computer
  Science, Jagiellonian University, \L ojasiewicza 6, 30-348 Krak\'ow,
  Poland} \thanks{The first- and second-named author were funded by
  the Jagiellonian University in Krak\'ow under ``Excellence Initiative
  – Research University''} \email{rafal.czyz@im.uj.edu.pl}

\author{P. H. Lundow}
\address{Department of Mathematics and Mathematical Statistics,
  Ume{\aa} University, SE-901 87 Ume{\aa}, Sweden}
\email{per-hakan.lundow@umu.se}


\title[Complex branches of a generalized Lambert $W$ function]{Complex
  branches of a generalized Lambert $W$ function arising from
  $p,q$-binomial coefficients}

\begin{abstract}
The $\psi(x)$-function, which solves the equation $x = \sinh(aw)e^w$ for $0<a<1$, has a natural connection to the renowned Lambert $W$ function and also physical relevance through its connection to the Lenz-Ising model of ferromagnetism. We give a detailed analysis of its complex branches and construct Riemann surfaces from these under various conditions of $a$, unveiling intriguing new links to the Lambert $W$ function.
\end{abstract}

\maketitle

\section{Introduction}
Special functions are the bridges between pure mathematics and real-world applications. A prime example is the Lambert $W$ function, the inverse function of $we^w$. Beyond its mathematical elegance, the Lambert $W$ function has captured the interest of mathematicians, physicists, and engineers due to its utility in solving real-world problems. Its applications range from modeling biological growth processes to optimizing electrical circuits, as evidenced in works like \cite{BhamidiSteeleZaman, Kozlov, TrefethenWeideman}. Furthermore, the function has inspired a variety of generalizations \cite{BariczMezo,ScottFreconGrotendorst,ScottMann}, each opening new avenues for exploration and innovation.

To provide some introductory background and context, we will briefly
describe how we arrived at the present article. The Lenz-Ising spin
model of ferromagnetism is one of the most celebrated theoretical physics models.
With $\pm 1$-spins located on the $d$-dimensional
integer lattice and nearest-neighbor interaction have been solved
completely only for the $1$-dimensional case ($d=1$) and partially
solved (without an external field) for $d=2$ in the 1940s. For $d=3$
very little is known rigorously, but good estimates of critical
parameters abound. For $d\ge 4$, the critical exponents, which
describe how relevant quantities behave around the critical
temperature $T_c$, are the same as for a mean-field model,
corresponding to spin interactions between the vertices of a complete
graph.

One important quantity, magnetization (an order parameter), is
the sum of the spins. As it happens, the magnetization
distribution for a complete graph is given by a $p,q$-binomial
distribution where $p=q=\exp(-2/T)$.  Also, for
$d\ge 5$, the magnetization distribution is extremely close to a
$p,q$-binomial distribution for finite $n$ and asymptotically the same
at $T_c$. However, it is far from clear how $p,q$ depends on the
temperature parameter $T$. Thus, the high-dimensional Lenz-Ising model
and the $p,q$-binomial coefficients are related, but it is not clear to
what extent~\cite{LundowRosengren}.

We remind the reader here that the $p,q$-binomial coefficients are
defined for $p\ne q$ as
\[\pqbin{n}{k} = q^{k(n-k)}\qbin{n}{k}{p/q},\]
where the right-hand side contains the standard $q$-binomial
coefficient.

For $0<p,q<1$, this distribution is symmetric and has its modes (peak
locations) at $k=(n/2)(1\pm a)$ where $0<a<1$, thus being unimodal
($a=1/2$) or bimodal ($a\ne 1/2$). The parameter $a$ corresponds to the Lenz-Ising model's normalized magnetization.
Choosing $0<a<1$ and letting $p=1+2y/n$ and $q=1+2z/n$, where $z$ can
be thought of as a distribution shape parameter and $y$ is the
solution to $\sinh(az)e^z=\sinh(ay)e^y$, now gives us a
$p,q$-binomial distribution with mode $a$ (when $n\to\infty$).

This means that we must solve $x=\sinh(ay)e^y$ for $y$ with
$x=\sinh(az)e^z$, i.e., we compute a branch of $\psi(a,x)$, an inverse
of $\sinh(az)e^z$. This corresponds to a mapping between different
branches of $\sinh(az)e^z$ which is a direct and natural
generalization of mapping between Lambert $W$ branches $W_0$ and
$W_{-1}$ obtained when $a=0$.  This process is described in detail in
\cite{LundowRosengren} and \cite{AhagCzyzLundow} where the latter also
gives some mathematical properties of the two real branches of $\phi$,
along with the transition function $\omega$ which maps between the two
real branches.

Note that the mathematical properties of the $p,q$-binomial
coefficients are interesting in their own right (see
\cite{AhagCzyzLundow} for references) even without their connection to
the Lenz-Ising model. However, we are not aware of any studies of
their properties for complex-valued $p,q$. This would then lead to
complex-valued probability amplitudes with the squared modulus as
probability density, as is fundamental in quantum mechanics. Before
such a study is undertaken, we will investigate the properties of
the complex branches of $\psi$, which is what we hope to achieve with
this article, just as in prior studies on the classic Lambert $W$
function~\cite{Beardon,Beardon2,CorlessGonnetHareJeffreyKnuth}.\\

Let $f(w)=\sinh(aw)e^w$ with $a$ treated as a fixed parameter
($0<a<1$) and $\psi$ as the inverse function.  We begin with the real
branches of $\psi$:

\begin{enumerate}\itemsep2mm
    \item $\psi_0$, an increasing and concave function defined as
    \[
    \psi_0:[x_a,+\infty) \to [\xi_a,+\infty),
    \]
    where $x_a =
    \frac{-a}{1+a}\left(\frac{1-a}{1+a}\right)^{\frac{1-a}{2a}}$ and
    $\xi_a = \frac{1}{2a}\ln \left(\frac{1-a}{1+a}\right)$.

    \item $\psi_{-1}$, a decreasing function defined as
    \[
    \psi_{-1}:[x_a,0) \to (-\infty,\xi_a].
    \]
\end{enumerate}

\medskip

In a complex generalization of $\psi$, giving a multi-valued inverse
of $f$, the behavior of its branches depends significantly on the
parameter $a$. In Section~\ref{sec:categ} we identify three categories
of $a$ which are then dealt with separately:

\medskip

\begin{enumerate}\itemsep2mm
\item $\frac {1+a}{1-a}\in \mathbb N$ (Section~\ref{sec:case(1)});

\item $\frac {1+a}{1-a}\notin \mathbb N$ and $a\in \mathbb Q$ (Section~\ref{sec:case(2)});

\item $a\notin \mathbb Q$ (Section~\ref{sec:case(3)}).

\end{enumerate}

\medskip

Classically, complex branches are connected to Riemann surfaces and
therefore, in Sections~\ref{sec:case(1)} and~\ref{sec:case(2)}, we aim
to construct these from the resultant complex branches, as has been
done earlier by Mez\H{o}~\cite{Mezo21} for a generalized version of
the Lambert function.

We will also treat the special cases $a\to 0^+$, $a\to 1^-$ and
$a=1/3$ separately.  For example, in the first case there is direct
relation to the classical Lambert $W$ function. More on the Lambert
$W$ function can be found in the highly recommended
monograph~\cite{MezoBook}.

\section{The Jacobian of $f$}

In this section, we investigate the behavior of $z=f(w) =
\sinh(aw)e^w$ where $z=x+iy$ and $w=\xi + i\eta$ with $x,y,\xi,\eta$
real.  We then get:
\begin{multline*}
z = x+iy = \sinh(aw)e^w = \frac{1}{2}\left(e^{(1+a)(\xi+i\eta)}-e^{(1-a)(\xi+i\eta)}\right)\\
  = \frac{1}{2}\Bigg(e^{(1+a)\xi}(\cos((1+a)\eta)+i\sin((1+a)\eta))\\ -e^{(1-a)\xi}(\cos((1-a)\eta)+i\sin((1-a)\eta))\Bigg),
\end{multline*}
resulting in
\[
\begin{aligned}\label{coordinates}
x &= \frac{1}{2}\left(e^{(1+a)\xi}\cos((1+a)\eta)-e^{(1-a)\xi}\cos((1-a)\eta)\right),\\
y &= \frac{1}{2}\left(e^{(1+a)\xi}\sin((1+a)\eta)-e^{(1-a)\xi}\sin((1-a)\eta)\right),\\
\end{aligned}
\]
with
\[
x^2+y^2 = \frac{1}{2}e^{2\xi}(\cosh(2\xi a)-\cos(2\eta a)).
\]
By differentiating these expressions, we can now calculate the Jacobian of
$f(w) = \sinh(aw)e^w$. We find that:
\[
\begin{aligned}
x_{\xi} &= \frac{1}{2}\left((1+a)e^{(1+a)\xi}\cos((1+a)\eta)-(1-a)e^{(1-a)\xi}\cos((1-a)\eta)\right),\\
x_{\eta} &= \frac{1}{2}\left(-(1+a)e^{(1+a)\xi}\sin((1+a)\eta)+(1-a)e^{(1-a)\xi}\sin((1-a)\eta)\right),\\
y_{\xi} &= \frac{1}{2}\left((1+a)e^{(1+a)\xi}\sin((1+a)\eta)-(1-a)e^{(1-a)\xi}\sin((1-a)\eta)\right) = -x_{\eta},\\
y_{\eta} &= \frac{1}{2}\left((1+a)e^{(1+a)\xi}\cos((1+a)\eta)-(1-a)e^{(1-a)\xi}\cos((1-a)\eta)\right) = x_{\xi},
\end{aligned}
\]
yielding
\begin{equation}\label{jacobian}
\begin{aligned}
J &= x_{\xi}y_{\eta}-x_{\eta}y_{\xi} = (x_{\xi})^2+(x_{\eta})^2\\
  &= \frac{1}{4}e^{2\xi}\left((1+a)^2e^{2a\xi}+(1-a)^2e^{-2a\xi}-2(1-a^2)\cos(2\eta a)\right)\geq 0.
\end{aligned}
\end{equation}
Solving the equation $J=0$, we get the following critical points:
\[
\begin{aligned}
(1+a)^2e^{2a\xi} &= (1-a)^2e^{-2a\xi}, \ \ \text{and}\\
\cos(2\eta a) &= 1,
\end{aligned}
\]
so that
\[
w_k=\xi+i\eta_k=\frac{1}{2a}\ln \left(\frac{1-a}{1+a}\right)+i\frac{k\pi}{a}, \ \ k\in \mathbb{Z},
\]
and then
\[
z_k=f(w_k)=\frac{-a}{1+a}\left(\frac{1-a}{1+a}\right)^{\frac{1-a}{2a}}(-1)^k\left(\cos\left(\frac{\pi k}{a}\right)+i\sin\left(\frac{\pi k}{a}\right)\right).
\]
We denote this set of critical points as $\textbf{CP}_a = \{z_k: k \in
\mathbb{Z}\}$. As we will see, $\textbf{CP}_a$ also represents the
branch points for $\psi$. By using the Implicit Function Theorem, we
can calculate the derivative of $\psi$ for each of its branches:
\begin{equation*}
  \psi'(z) = \frac{1}{e^{\psi(z)}\left(a\cosh(a\psi(z))+\sinh(a\psi(z))\right)}.
\end{equation*}
Note that $\psi'$ is not defined at the critical points $z_k$.


\section{The Principal Complex Branch $\psi_0$}\label{sec:principal}

The principal complex branch $\psi_0$ has a consistent appearance for
every parameter $a\in (0,1)$. Our goal is to generate a complex
extension of the real branch $\psi_0$ such that the interval
$(x_a,+\infty)$ is included within the domain of definition $\Omega_0$
of $\psi_0$.

The first step involves locating the points at which
$y=\operatorname{Im}f(z)=0$. From (\ref{coordinates}), we have $y=0$
if and only if $\eta=0$ or $\eta=\frac {k\pi}{1-a}$, for $k\in \mathbb
Z$ (if $\frac {1+a}{1-a}\in \mathbb N$), or if
\[
e^{(1+a)\xi}\sin((1+a)\eta)-e^{(1-a)\xi}\sin((1-a)\eta)=0.
\]
Let us designate the solutions to this equation as
\begin{equation}\label{3.1}
\Xi(\eta)=\xi=\frac {1}{2a}\ln\left(\frac {\sin((1-a)\eta)}{\sin((1+a)\eta)}\right).
\end{equation}

\begin{figure}[!ht]
    \centering
    \includegraphics[width=0.85\textwidth]{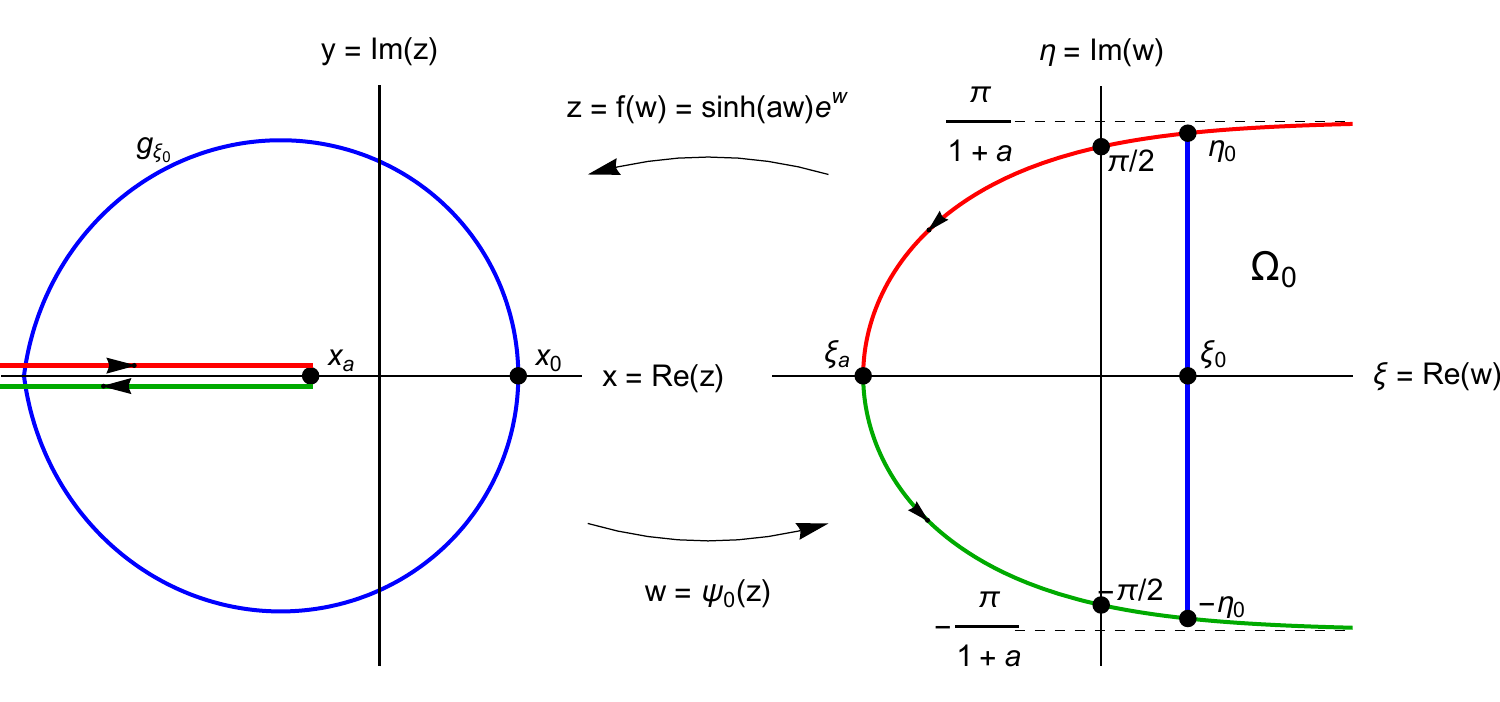}
    \caption{Construction of the principal branch $\psi_0$ for any
      parameter $a\in (0,1)$. The domain $\mathbb C\setminus
      (-\infty,x_a\rbrack$ (Left) is mapped by $\psi_0$ into
      $\Omega_0$ (Right). The red and green segments
      $(-\infty,x_a\rbrack$ on the left correspond to the red and
      green lines on the right, as given by the formula
      (\ref{3.1}). The arrows indicate the direction of movement along
      the curves. The vertical segment $\lbrack -\eta_0,\eta_0\rbrack$
      on the right is mapped by $f(w)=\sinh(aw)e^w$ onto the blue
      curve $g_{\xi_0}$ on the left, as described by the
      parametrization (\ref{3.3}).}
    \label{IMG:Nr1}
\end{figure}

The function $\Xi(\eta)$ is well-defined when $\frac
{\sin((1-a)\eta)}{\sin((1+a)\eta)}>0$. Let us denote $I_0=(-\frac
{\pi}{1+a},\frac {\pi}{1+a})$ and
\[
\Omega_0=\left\{(\xi,\eta): \xi>\Xi(\eta)=\frac {1}{2a}\ln\left(\frac {\sin((1-a)\eta)}{\sin((1+a)\eta)}\right), \eta\in I_0\right\}.
\]
We propose that
\[
\psi_0:\mathbb C\setminus(-\infty,x_a]\to \Omega_0
\]
is a conformal bijection. The proof of this involves observing that
the function $f(w)=\sinh(aw)e^w$ maps the upper half of the graph of
$\Xi(\eta)$ for $\eta\in \lbrack 0,\tfrac {\pi}{1+a})$, namely the set
\[
\left\{\left(\frac {1}{2a}\ln\left(\frac {\sin((1-a)\eta)}{\sin((1+a)\eta)}\right),\eta\right): \eta\in \left\lbrack 0,\frac{\pi}{1+a}\right)\right\},
\]
to the set $(-\infty,x_a\rbrack$. This can be seen from relation
(\ref{coordinates}), where we have $y=0$ and
\begin{equation}\label{x}
  x=\frac{1}{2}e^{\xi(1-a)}\frac {-\sin(2\eta a)}{\sin(\eta(1+a))}=\frac {-\sin(2\eta a)}{2\sin(\eta(1+a))}\left(\frac {\sin((1-a)\eta)}{\sin((1+a)\eta)} \right)^{\frac {1-a}{2a}}<0.
\end{equation}
The maximum value is achieved at $\eta=0$ with a value of $x_a$. This
observation can be mirrored for the lower/negative half of the graph
of $\Xi$. To validate the injectivity of $f(w)=\sinh(aw)e^w$ within
$\Omega_0$, remember that the Jacobian of $f$ is positive in
$\Omega_0$. It is worth noting that for a fixed $\xi_0$, the
horizontal interval $\left(\{\xi_0\}\times \mathbb R\right)\cap
\Omega_0$ is mapped by function $f$ into a simple curve $g_{\xi_0}$
(distinct for different $\xi$):

\begin{equation}\label{3.3}
\begin{aligned}
&x=\alpha(\xi_0)\cos(\eta(1+a))-\beta(\xi_0)\cos(\eta(1-a)),\\
&y=\alpha(\xi_0)\sin(\eta(1+a))-\beta(\xi_0)\sin(\eta(1-a)),\\
\end{aligned}
\end{equation}
where $\alpha(\xi)=\frac{1}{2}e^{\xi(1+a)}$ and
$\beta(\xi)=\frac{1}{2}e^{\xi(1-a)}$, see Fig.~\ref{IMG:Nr1}.

\section{Categorizing the Complex Branches}\label{sec:categ}

We define the following sets: $\mathbb C_{+}=\{z\in \mathbb C:
\operatorname{Im} z>0\}$, $\mathbb C_{-}=\{z\in \mathbb C:
\operatorname{Im} z<0\}$, and
\[
\Omega_k=\Omega_0+\frac {ik\pi}{a}, \ \ k\in \mathbb Z.
\]
The function $f$ exhibits a particular property: given $\xi_k=\xi_0$
and $\eta_k=\eta_0+\frac{k}{a}\pi $ with $(\xi_0,\eta_0)=w_0\in
\Omega_0$, we have
\begin{equation*}
f(w_k)=f(\xi_k,\eta_k)=f(w_0)e^{\frac {a+1}{a}k\pi i}.
\end{equation*}
Consequently, our investigation can focus on the behavior of $f$ in
the strip $\mathbb R\times [0,\frac{\pi}{a})$. If $a=\frac {p}{q}\in
  \mathbb Q$, then $f$ is a periodic function with a period defined as
\[
T_a=\left\{
\begin{array}{ll}
  q\pi i, & p+q\in 2\mathbb N \\
  2q\pi i, & p+q\notin 2\mathbb N.
\end{array}
\right.
\]
The behavior of other branches, however, is more intricate and depends
on the properties of the function
\[
\Xi(\eta)=\frac {1}{2a}\ln\left(\frac {\sin((1-a)\eta)}{\sin((1+a)\eta)}\right).
\]
We divide the $\eta$-axis into disjoint intervals, distinguished by
points at $\frac {k\pi}{1-a}$ and $\frac {k\pi}{1+a}$, $k\in \mathbb
Z$. Here, both sine functions equal zero. Despite the tempting
possibility of only considering intervals $I_k=(\frac {k\pi}{1+a},
\frac {k\pi}{1-a})$, these intervals will eventually intersect for
large $k$. We note that $\Xi(\eta)$ does not possess double zeros,
implying that $\frac {\sin((1-a)\eta)}{\sin((1+a)\eta)}$ alternates in
sign across intervals. The function remains even for all cases. We
denote $I_k=(a_k,b_k)$, where $a_k,b_k\in \left\{\frac {k\pi}{1-a},
\frac {k\pi}{1+a}, k\in \mathbb Z\right\}$, as the intervals where
$\Xi(\eta)$ is well defined. Thus,
\[
\dom(\Xi)=\bigcup_{k}I_k.
\]
Examining the behavior of $\Xi$ more carefully, we find extreme points
of $\Xi$ when $\Xi'(\eta)=0$. This holds if and only if
\[
(1-a)\cos(\eta(1-a))\sin(\eta(1+a))=(1+a)\cos(\eta(1+a))\sin(\eta(1-a)),
\]
yielding
\[
\sin(2\eta a)=a\sin(2\eta).
\]
Let $\eta_k'$ represent the solutions to the above equation. If $\frac
{1+a}{1-a}\in \mathbb N$, then the only solutions are $\eta_k'=\frac
{k\pi}{1-a}$, where $\Xi$ attains its minimum. For the case when
$\frac {1+a}{1-a}\notin \mathbb N$, the situation is similar. However,
there may exist points at $\eta_k'=\frac {k\pi}{1-a}$ where
$\Xi'(\eta_k')>0$, thereby not reaching an extremum.

Let us consider the points where $\Xi(\eta)=0$. This is equivalent to
$\sin(\eta(1-a))=\sin(\eta(1+a))$, leading to $\sin(\eta a)\cos
\eta=0$, and hence $\eta=\frac{\pi k}{a}$ or
$\eta=\frac{\pi}{2}+k\pi$, for $k\in \mathbb Z$.

As for the behavior at infinity, let us consider the interval $I_k$,
where $\Xi$ is well defined. The boundary points of this interval
correspond to zeros of either the numerator, $\sin((1-a)\eta)$, or the
denominator, $\sin((1+a)\eta)$. In this scenario, we may observe two
distinct cases:

\begin{enumerate}\itemsep2mm
\item All zeros of the numerator $\sin((1-a)\eta)$ coincide with zeros
  of the denominator $\sin((1+a)\eta)$, i.e., $\frac{1+a}{1-a}\in
  \mathbb N$. In such a case, as $\eta$ approaches the boundary of
  $I_k$, $\Xi(\eta)$ tends to $+\infty$.

\item Not all zeros of the numerator $\sin((1-a)\eta)$ coincide with
  zeros of the denominator $\sin((1+a)\eta)$, i.e.,
  $\frac{1+a}{1-a}\notin \mathbb N$. In this situation, as $\eta$
  approaches the boundary of $I_k$, $\Xi(\eta)$ tends to $+\infty$
  (when the boundary point is a zero of the denominator) or tends to
  $-\infty$ (when the boundary point is a zero of the numerator but
  not the denominator).

\end{enumerate}

Based on the above considerations, the behavior of the branches of
$\psi$ depends on the parameter $a$ and can be categorized into three
scenarios, when

\begin{enumerate}\itemsep2mm
\item $\frac {1+a}{1-a}\in \mathbb N$;

\item $\frac {1+a}{1-a}\notin \mathbb N$, but $a\in \mathbb Q$;

\item $a\notin \mathbb Q$.
\end{enumerate}
\vfill

\section{The Case $\frac {1+a}{1-a}\in \mathbb N$}\label{sec:case(1)}

We start with the special case $a=\frac{1}{2}$ and then the general
case $\frac {1+a}{1-a}\in \mathbb N$.\\

\noindent\emph{Case $a=\frac{1}{2}$.}\\

\begin{figure}[!ht]
    \centering
    \includegraphics[width=0.80\textwidth]{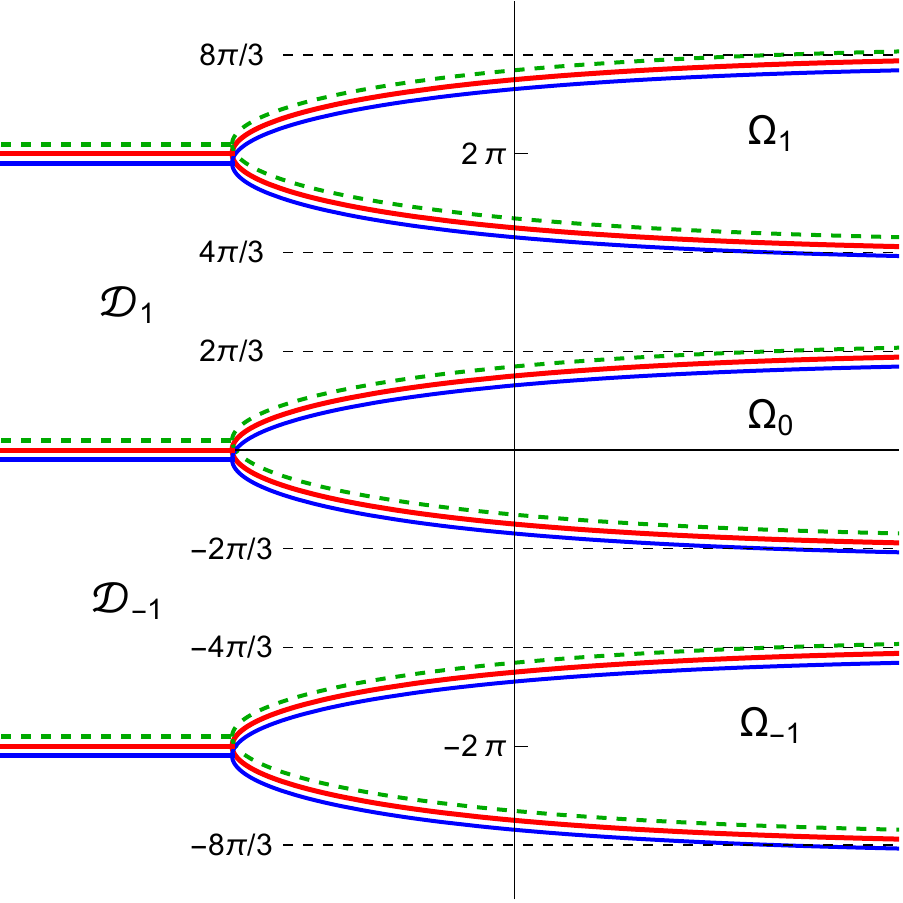}
    \caption{Division of the complex plane $\mathbb C_w$ for the
      parameter $a=\frac{1}{2}$ onto regions where $f(w)=\sinh(aw)e^w$
      is injective (the co-domains of the complex branches of
      $\psi_k$, $\tilde{\psi}_k$, as per (\ref{5.2}),
      (\ref{5.3})). The plot of $\Xi(\eta)=\frac {1}{2a}\ln\left(\frac
      {\sin((1-a)\eta)}{\sin((1+a)\eta)}\right)$ is shown in red. For
      the extension of the complex branches of $\psi$, refer to
      (\ref{5.4}): dashed green lines indicate parts of the boundary
      that do not belong to the corresponding co-domain, while solid
      blue lines indicate parts of the boundary that do belong to the
      corresponding co-domain.}
    \label{IMG:Nr2}
\end{figure}

If we define $\xi_k=\xi_0$ and $\eta_k=\eta_0+2k\pi i$, where
$(\xi_0,\eta_0)=w_0\in \Omega_0$, then we have
\begin{equation}\label{branches}
  f(w_k)=f(\xi_k,\eta_k)=(-1)^kf(w_0).
\end{equation}
Additionally, $f$ is a periodic function exhibiting a period of
$T_{1/2}=4\pi i$. The critical points for this case are
represented by $\textbf{CP}_{1/2}=\{x_{1/2},
-x_{1/2}\}$, where $x_{1/2}=-\frac {\sqrt 3}{9}$. The
following objective is to partition the $w$-plane into domains within
which $f(w)=\sinh(aw)e^w$ retains injective properties. With this goal
in mind, define
\[
\Omega_k=\Omega_0+2k\pi i, \ \ k\in \mathbb Z;
\]
and
\[
D_k=\mathbb R\times(2(k-1)\pi,2k\pi)\setminus \overline{(\Omega_{k-1}\cup\Omega_k)}.
\]
See Figure~\ref{IMG:Nr2} for an illustration of these sets.

It is important to note that all $\Omega_k$ and $D_k$ remain
open. Applying (\ref{branches}), we can deduce that $\psi$ has
countably many branches $\psi_k$ and $\widetilde\psi_k$ conformally
bijecting onto their corresponding domains. In specific terms,

\begin{enumerate}\itemsep2mm
\item $\psi$ has countably many branches $\psi_k$
  \begin{equation}\label{5.2}
    \begin{aligned}
      &\psi_{2k}:\mathbb C\setminus(-\infty,x_a]\to \Omega_{2k}, \quad\text{for } k\in \mathbb Z;\\
&\psi_{2k+1}:\mathbb C\setminus[-x_a,+\infty)\to \Omega_{2k+1}, \quad\text{for } k\in \mathbb Z;\\
    \end{aligned}
  \end{equation}

\item $\psi$ has countably many branches $\widetilde\psi_k$
  \begin{equation}\label{5.3}
    \begin{aligned}
      &\widetilde\psi_{k}:\mathbb C_{+}\to D_{k}, \quad\text{for }  k\in (2\mathbb N\setminus\{0\})\cup(-2\mathbb N-1);\\
      &\widetilde\psi_{k}:\mathbb C_{-}\to D_{k}, \quad\text{for } k\in (2\mathbb N+1)\cup(-2\mathbb N\setminus \{0\});\\
    \end{aligned}
  \end{equation}
\end{enumerate}

Furthermore, the domain of definition of the functions can be extended as follows:
\begin{equation}\label{5.4}
  \begin{aligned}
    &(1) \ \psi_k \ \text{ can be extended to encompass the entirety of } \ \mathbb C.\\
    &(2) \ \tilde \psi_k \ \text{ can be extended to }  \ \mathbb C_{+}\cup(-\infty,0) \text{ for } k\in (2\mathbb N\setminus\{0\})\cup(-2\mathbb N-1).\\
    &(3) \ \tilde \psi_k \ \text { can be extended to } \ \mathbb C_{-}\cup(0,\infty) \text{ for }  k\in (2\mathbb N+1)\cup(-2\mathbb N\setminus \{0\}).
  \end{aligned}
\end{equation}

It is worth noting that the real branches $\psi_0$ and $\psi_{-1}$ are
incorporated into the corresponding complex branches $\psi_0$ and
$\widetilde\psi_{-1}$ respectively. Examining the branch structure
reveals the presence of three branch points:
\[
\textbf{BP}_{1/2}=\textbf{CP}_{1/2}\cup\{0\}=\left\{-\frac {\sqrt 3}{9},\frac {\sqrt 3}{9},0\right\}.
\]

\begin{figure}[!ht]
    \centering
    \includegraphics[width=0.60\textwidth]{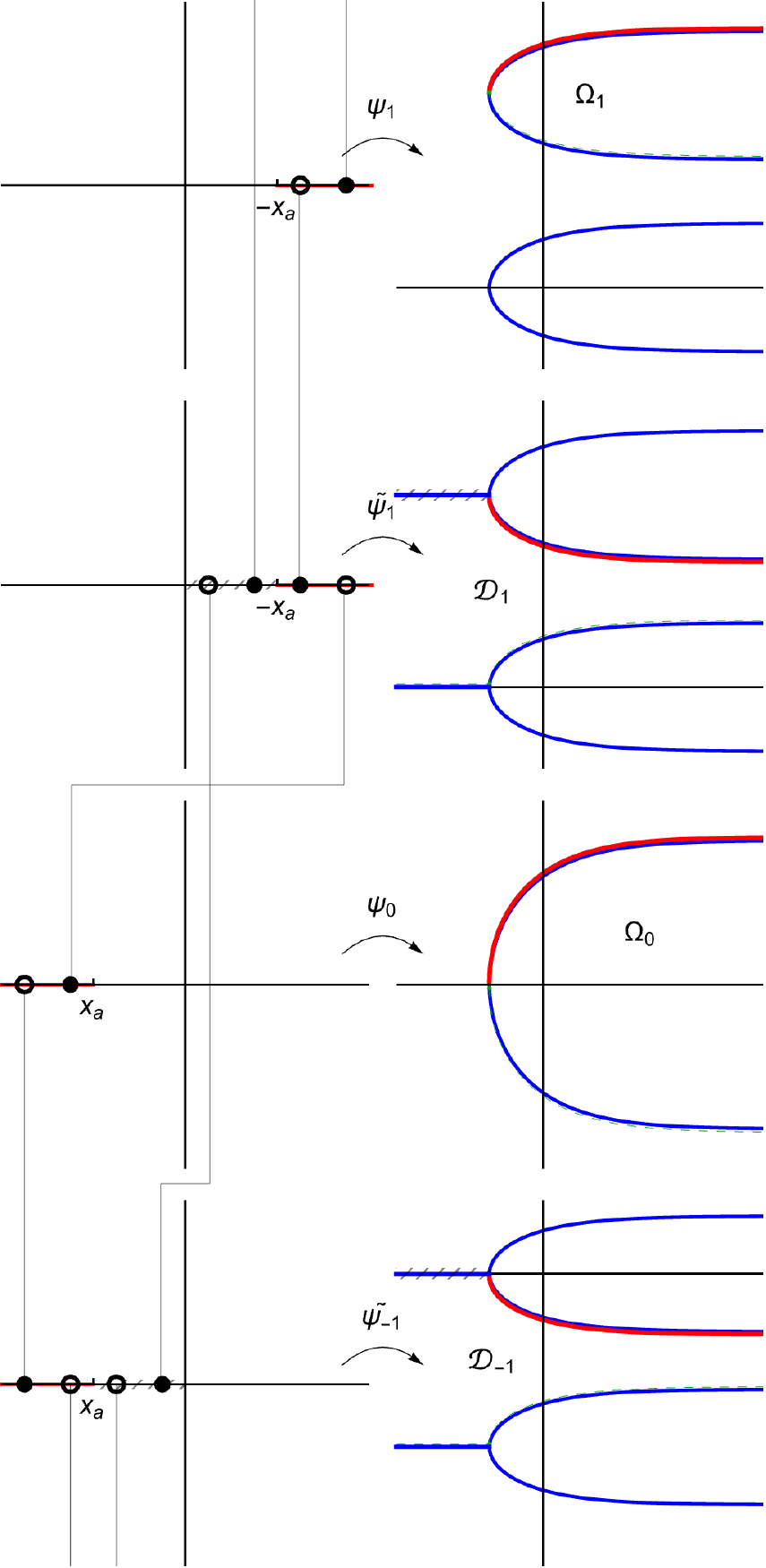}
    \caption{Construction of the Riemann surface for the complex
      branches of $\psi$, with the parameter $a=\frac{1}{2}$. The red
      segments (Left) correspond to the red curves (Right), while the
      ``crossed out'' black segments (Left) correspond to the ``crossed
      out'' blue segments (Right). In the left panel, during the
      construction of the Riemann surface, red segments and ``crossed
      out'' segments are meant to be ``glued'' together, respectively. A
      filled dot indicates that the corresponding cut is closed in the
      appropriate leaf of the Riemann surface, while an empty dot
      indicates that the cut is open.}
    \label{IMG:Nr3}
\end{figure}

Each of these branch points involves a specific combination of the
$\psi$ and $\widetilde \psi$ branches. Further details regarding this
interaction can be inferred from Figure ~\ref{IMG:Nr3}, which highlights
how intervals under various branches divide the entire $w$-plane,
where the ``upper'' boundary is added to the domain below.

Now let us focus on the construction of the Riemann surface associated
with $\psi$. The first step is to sever the domains of $\tilde
\psi_{-1}$ and $\tilde \psi_{1}$ along the intervals $\lbrack x_{1/2},
0\rbrack$ and $\lbrack 0,-x_{1/2}\rbrack$, respectively. Following
this, the cuts are pieced together such that the domain of $\psi_{-1}$
is closed at the interval while the domain of $\tilde \psi_1$ remains
open. This procedure is followed by a series of cutting and gluing
operations in a particular sequence. We refer to Figure~\ref{IMG:Nr3}
for a more detailed illustration.

Investigating the monodromy group around the branch points reveals the
group to be the infinite cyclic group, $\mathbb Z$, for each branch
point. This is inferred from observing the lifts of positively
oriented, closed curves with maximal modulus less than $|x_{1/2}|$ and
infinite winding numbers around the branch points. Finally, the chosen
branches are observed to satisfy the Counter Clock Continuity (CCC)
rule around each branch point (see Fig.~\ref{IMG:Nr3}).

\medskip

\noindent\emph{The general case.}\\

Next, we continue with the general case when $\frac {1+a}{1-a}\in
\mathbb N$. This implies that $a=\frac{n-1}{n+1}$ for a specific $n\in
\mathbb N$. Under these conditions, the function $f$ exhibits a
characteristic property: if $\xi_k=\xi_0$ and
$\eta_k=\eta_0+\frac{k}{a}\pi i$, where $(\xi_0,\eta_0)=w_0\in
\Omega_0$, it can be expressed as:
\begin{equation}\label{branches2}
  f(w_k)=f(\xi_k,\eta_k)=f(w_0)e^{\frac {a+1}{a}k\pi i}=f(w_0)e^{\frac{2n}{n-1}k\pi i}.
\end{equation}

This equation suggests that the behavior of $f$ within the strip
$\mathbb{R}\times\lbrack 0,\frac{\pi(n+1)}{n-1})$ is sufficient for
the investigation. The function $f$ also displays periodic behavior
with a period of $T_a=(n+1)\pi i$. The set of critical points can be
defined as:
\[
\textbf{CP}_a=\left\{z_k=x_a(-1)^k\left(\cos\left(\tfrac {\pi k}{a}\right)+i\sin\left(\tfrac {\pi k}{a}\right)\right): k=0,...,n-2\right\}.
\]
For the analysis, let us establish the following definitions:
\begin{align*}
  \Omega_k &= \Omega_0+\tfrac{k}{a}\pi i, \quad k\in \mathbb Z, \\
  D_k &= \left(\mathbb R\times\left(\tfrac{k-1}{a}\pi,\tfrac{k}{a}\pi\right)\right)\setminus \overline{(\Omega_{k-1}\cup\Omega_k)}.
\end{align*}
Both $\Omega_k$ and $D_k$ are open sets. Their partitioning of the
$w$-plane is reminiscent of the case when $a=\frac{1}{2}$, where one
whole period is covered by including $\Omega_0, \cdots,\Omega_{n-1}$
and $D_1,\cdots,D_{n-1}$. Now, consider the line segments:
\[
(-\infty, z_k\rbrack=\{e^{\frac {a+1}{a}k\pi i}x:\,  x\in (-\infty,x_a\rbrack\}, \quad k=0,\dots,n-2,
\]
and for $1\leq k\leq n-1$ let $\mathcal A_k$ denote the angle defined
by $(-\infty, z_{k-1}\rbrack$ and $(-\infty, z_{k}\rbrack$. The case
$a=\frac{1}{2}$ has only two critical points $x_{1/2}$ and $-x_{1/2}$,
where $\mathcal A_1=\mathbb C_-$ and $\mathcal A_2=\mathbb C_+$.

Leveraging arguments similar to those used for the principal branch
and the case $a=\frac{1}{2}$, and incorporating equation
(\ref{branches2}), we can prove that:
\begin{enumerate}\itemsep2mm
\item $\psi$ comprises an infinite set of branches $\psi_k$, each of
  which is a conformal bijection:
  \[
  \psi_{k}:\mathbb C\setminus (-\infty,z_k\rbrack\to \Omega_{k},
  \quad \text {for } k\in \mathbb Z,
  \]

  \item $\psi$ also has an infinite set of branches
    $\widetilde\psi_k$, each of which is also a conformal bijection:
    \[
    \widetilde\psi_{k}:\mathcal A_k\to D_{k}, \quad\text{for } k\in \mathbb Z.
    \]
\end{enumerate}

In order to segment the entire $w$-plane, the boundary must be
adjusted as follows: The ``upper'' boundary is appended to the domain
below in a manner similar to the case $a=\frac{1}{2}$. As a result,
the domain of definition of the function can be expanded:
\begin{enumerate}\itemsep2mm
\item $\psi_k$ can be extended to the whole $\mathbb C$,
\item $\tilde \psi_k$ can be extended to $\mathcal A_k\cup \{tz_k: t\in(-\infty,0)\}$.
\end{enumerate}
It is noteworthy that the real branches $\psi_0$ and $\psi_{-1}$ are
subsumed within the corresponding complex branches $\psi_0$ and
$\widetilde\psi_{-1}$, respectively. The structure of these branches
reveals that the branch points correspond to:
\[
\textbf{BP}_{a}=\textbf{CP}_{a}\cup\{0\}.
\]

The process to construct the Riemann surface associated with $\psi$
and the description of the monodromy group parallels the approach for
the case of $a=\frac{1}{2}$. Due to this similarity, we omit a
detailed discussion. It is essential to note, however, that while the
overall construction procedure remains the same, the number of ``gluing
parts'' will vary.


\section{The Case $\frac {1+a}{1-a}\notin \mathbb N$, $a\in \mathbb Q$}\label{sec:case(2)}
As in the previous section we will start with a specific value of $a$,
and then afterwards proceed to the general case.\\

\noindent\emph{Case $a=\frac{1}{4}$.}\\

Let $a=\frac{1}{4}$, and consider the $f(w) = \sinh(\frac{1}{4}w)e^w$
which is periodic with a period of $8\pi i$. This function has
critical points $\textbf{CP}_{1/4}=\{x_{1/4},
-x_{1/4}\}$, where $x_{1/4}=-\frac{3\sqrt{15}}{125}$.

For this case, we need to partition the $w$-plane into domains where
$f(w)=\sinh(aw)e^w$ is injective. This scenario is different from the
case $a=\frac{1}{2}$, because the graph of the function $\Xi$ contains
``cubic'', denoted by $\Gamma_k$ for $k=1,2,5,6$, in addition to
the ``parabolic'' curves, denoted by $\Gamma_k$ for $k=0,3,4,7$.

Let us  divide the domain of $\Xi$ on $\lbrack 0,8\pi)$ into intervals $J_k$:

\begin{align*}
  J_0 & = \Big(0, \frac{4}{5}\pi\Big),  & J_1 & = \Big(\frac{4}{3}\pi, \frac{8}{5}\pi\Big),  & J_2 & = \Big(\frac{12}{5}\pi, \frac{8}{3}\pi\Big),  & J_3 & = \Big(\frac{16}{5}\pi, 4\pi\Big), \\
  J_4 & = \Big(4\pi, \frac{24}{5}\pi\Big), & J_5 & = \Big(\frac{16}{3}\pi, \frac{28}{5}\pi\Big), & J_6 & = \Big(\frac{32}{5}\pi, \frac{20}{3}\pi\Big), & J_7 & = \Big(\frac{36}{5}\pi, 8\pi\Big).
\end{align*}

\begin{figure}[!ht]
    \centering
    \includegraphics[width=0.75\textwidth]{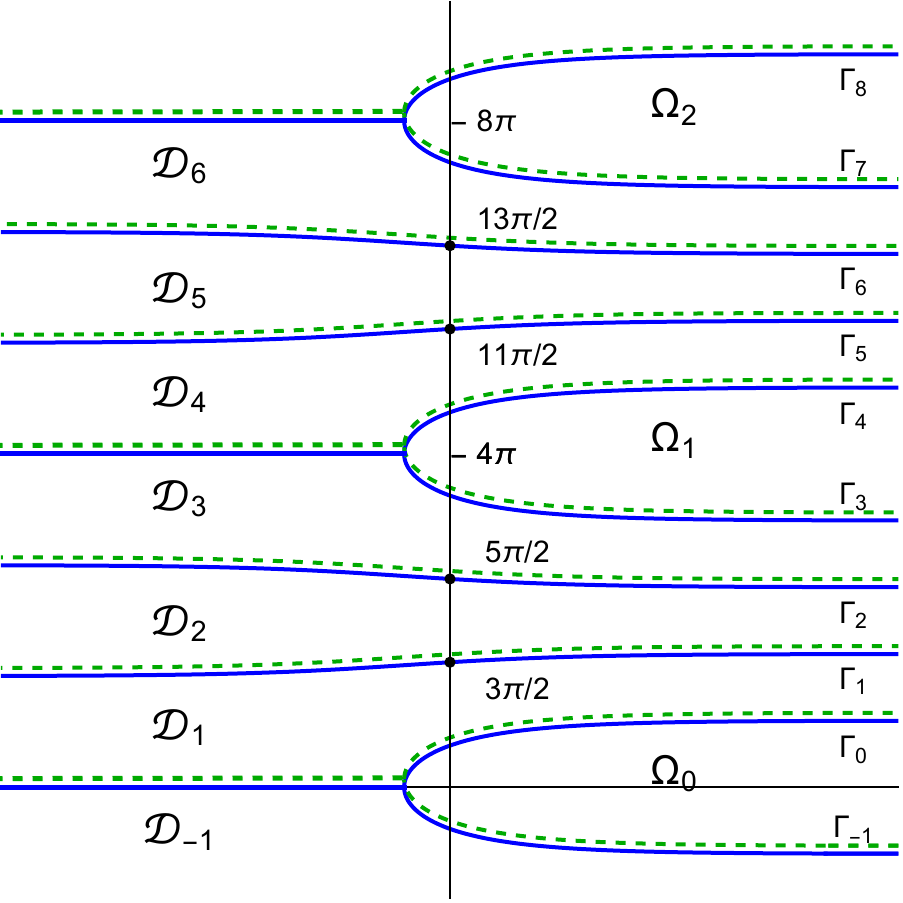}
    \caption{This figure illustrates the division of the complex plane
      $\mathbb C_w$ for the parameter $a=\frac{1}{4}$ onto regions
      where the function $f(w)=\sinh(aw)e^w$ is injective. These
      regions correspond to the co-domains of the complex branches of
      $\psi_k$ and $\tilde{\psi}_k$ (see (\ref{6.1}),
      (\ref{6.2})). The blue curve represents the plot of
      $\Xi(\eta)=\frac {1}{2a}\ln\left(\frac
      {\sin((1-a)\eta)}{\sin((1+a)\eta)}\right)$, and the curves
      $\Gamma_j$ for $j=-1,\cdots, 8$. Regarding the extension of the
      complex branches of $\psi$, refer to (\ref{6.3}). The dashed
      green lines indicate that the corresponding parts of the
      boundary do not belong to the corresponding co-domain, while the
      solid blue lines show that they do.} \label{IMG:Nr4}
\end{figure}

Here, the curve $\Gamma_k$ is defined as the graph of $\Xi$ on the
interval $J_k$ for $k=1,2,5,6$ and the graph of $\Xi$ on the interval
$J_k$ plus the interval $(-\infty, x_{1/4}\rbrack$ for
$k=0,3,4,7$. This is illustrated in Figure~\ref{IMG:Nr4}.

We define $\Omega_k$ as $\Omega_0+4k\pi$ for $k\in \mathbb{Z}$, and
$D_k$ as the regions between the curves $\Gamma_{k-1}$ and
$\Gamma_k$. Note that all $\Omega_k$ and $D_k$ are open. The function
$\psi$ has a countable number of branches $\psi_k$ and
$\widetilde\psi_k$, defined as follows:

\begin{enumerate}\itemsep2mm
\item  Branches of $\psi$:
  \begin{equation}\label{6.1}
    \begin{aligned}
      &\psi_{2k}: \mathbb{C}\setminus(-\infty,x_a\rbrack \rightarrow
      \Omega_{2k}, \quad\text{for } k\in \mathbb{Z},\\ &\psi_{2k+1}:
      \mathbb{C}\setminus\lbrack -x_a,+\infty) \rightarrow
      \Omega_{2k+1}, \quad\text{for } k\in \mathbb{Z},
    \end{aligned}
  \end{equation}

\item Branches of $\widetilde\psi$:
  \begin{equation}\label{6.2}
    \begin{aligned}
      &\widetilde\psi_{k}: \mathbb{C}_{+}  \rightarrow D_{k}, \quad\text{for }  k\in (2\mathbb{N}\setminus\{0\})\cup(-2\mathbb{N}-1),\\
      &\widetilde\psi_{k}: \mathbb{C}_{-}  \rightarrow D_{k},  \quad\text{for }  k\in (2\mathbb{N}+1)\cup(-2\mathbb{N}\setminus \{0\}),
    \end{aligned}
  \end{equation}
\end{enumerate}
where each $\psi_k$ and $\widetilde \psi_k$ forms a conformal bijection.

Now let us discuss the boundary behavior. According to
Equation~(\ref{x}), $f$ maps $\Gamma_k$ into $(0,+\infty)$ for
$k=1,3,4,6$ and into $(-\infty,0)$ for $k=0,2,5,7$. The entire
$w$-plane can be divided by the curves $\Gamma_k$ such that the
``upper'' boundary is included in the domain below. Therefore, the
domain of definition of the function

\begin{equation}\label{6.3}
  \begin{aligned}
    &(1) \, \psi_k  \text { can be extended to the whole }  \mathbb{C}\\
    &(2) \, \widetilde\psi_k \text { can be extended to }
    \mathbb{C}_{+}\cup(-\infty,0) \text{ for } k=2,4,6 \\
    &(3) \, \widetilde\psi_k \text { can be extended to }
    \mathbb{C}_{-}\cup(0,\infty) \text { for } k=1,3,5.\\
  \end{aligned}
\end{equation}

The real branches $\psi_0$ and $\psi_{-1}$ are included in the
corresponding complex branches $\psi_0$ and $\widetilde\psi_{-1}$
respectively.  The branch structure reveals that there are three
branch points denoted by
\[
\textbf{BP}_{1/4}=\textbf{CP}_{1/4}\cup\{0\}=\left\{-\frac {3\sqrt {15}}{125},\frac {3\sqrt {15}}{125},0\right\}.
\]

The branch point $x_{1/4}$ involves all $\psi_1$, $\tilde \psi_1$ and
$\tilde \psi_{-1}$ branches. The other branch point $0$ encompasses
all branches $\tilde \psi_k$, for $k=1,\dots,6$. Finally, the branch
point $-x_{1/4}$ encompasses all $\psi_1$, $\tilde \psi_3$, and
$\tilde \psi_{4}$ branches.

We now describe the Riemann surface associated with $\psi$. This
construction closely follows the case $a=\frac{1}{2}$, with the
exception that we now have three branches $\tilde \psi_k$, for
$k=1,2,3$ between $\psi_0$ and $\psi_1$ when $a=\frac{1}{4}$.
Initially, we take the interval $(-\infty,x_{1/4}\rbrack$ from
$\dom(\psi_0)$ and attach it to the interval $\lbrack
-x_{1/4},\infty)$ from $\dom(\tilde \psi_1)$, ensuring that
$\dom(\psi_0)$ is closed and $\dom(\tilde \psi_1)$ is open at the cut.

Subsequently, we make cuts along the interval $(0,+\infty)$ in the
domain of $\tilde \psi_{1}$ and along the interval $(-\infty,0)$ in
the domain of $\tilde \psi_{2}$, and then merge these cuts such that
$\dom(\psi_{1})$ is closed at the interval while $\dom(\tilde \psi_2)$
is open. Similarly, we attach the interval $(-\infty,0)$ of the domain
$\tilde \psi_{2}$ with the interval $(0,\infty)$ of the domain $\tilde
\psi_{3}$, where the first cut is closed while the second is open.

Eventually, we cut the domain of $\tilde\psi_{3}$ along the line
$\lbrack -x_{1/4},\infty)$ and glue it to the interval $\lbrack
-x_{1/4},\infty)$ of the domain of $\psi_1$ in such a way that
$\dom(\tilde\psi_{3})$ is closed and $\dom(\psi_1)$ is open at the
cut. The remaining closed cut of $\dom(\tilde \psi_{3})$ at the
interval $(0,-x_{1/4})$ is connected with an open cut of $\dom(\tilde
\psi_4)$ at the interval $(x_{1/4},0)$. This process is continued with
domains of $\tilde \psi_4$, $\tilde \psi_5$ and $\tilde\psi_6$. Our
choice of branches satisfy the Counter Clock Continuity (CCC) rule
around the branch point.

Regarding the monodromy group around the branch points, consider a
positively oriented, closed curve of maximum modulus smaller than
$|x_{1/4}|$, which winds infinitely many times around the branch point
$0$. The lift of this curve on the Riemann surface passes through the
sheets, excluding all $\psi_k$. Consequently, the monodromy group is
the infinite cyclic group, denoted as $\mathbb Z$.  At the points
$x_{1/4}$ and $-x_{1/4}$, the monodromy behaves similarly, but in this
case, $\psi_k$ with $k\in 2\mathbb Z$ and $k\in 2\mathbb Z+1$
respectively are also involved, thus reinforcing that the monodromy
group is indeed $\mathbb Z$.\\

\noindent\emph{The general case.}\\

Now we continue with the general case when $\frac {1+a}{1-a}\notin
\mathbb N$, $a\in \mathbb Q$. Here, $f$ represents a periodic function
with the period $T_a$ defined as:
\[
T_a=\begin{cases}
q\pi i, & \text{if } p+q \in 2\mathbb N, \\
2q\pi i, & \text{if } p+q \notin 2\mathbb N.
\end{cases}
\]
The set of critical points is given by:
\[
\textbf{CP}_a=\left\{z_k=x_a(-1)^k\left(\cos\left(\tfrac{\pi k}{a}\right)+i\sin\left(\tfrac{\pi k}{a}\right)\right) \mid k\in \mathbb Z\right\}.
\]
This situation presents a higher level of complexity, and a general
description is not readily available. The analysis bears some
similarity to the case $\frac{1+a}{1-a}\in \mathbb N$ and
$a=\frac{1}{4}$. However, the partitioning of the $w$-plane by the
curves $\Gamma_k$ described in the case $a=\frac{1}{4}$ is somewhat
different. The placement of ``parabolic'' and ``cubic'' curves $\Gamma_k$
depends on the parameter $a$.

We define $\Omega_k$ as follows:
\[
\Omega_k=\Omega_0+\frac{k}{a}\pi i, \quad k\in \mathbb Z,
\]
and $D_k$ as the regions between curves $\Gamma_{k-1}$ and
$\Gamma_k$. Please note that all $\Omega_k$ and $D_k$ are open. Let
$\mathcal{A}_k$ denote the angle between intervals $(-\infty,
z_{k-1}\rbrack$ and $(-\infty, z_{k}\rbrack$ for $k\in \mathbb
Z$. Define $\mathbb C_{k+}$ and $\mathbb C_{k-}$ as the upper and
lower half-planes with the boundary line passing through $0$ and $z_k$
(for $k=0$, we have the standard $\mathbb C_+$ and $\mathbb C_-$). In
this case, we identify three distinct types of branches: $\psi_k$,
$\widetilde{\psi}_k$, and $\widehat{\psi}_k^\pm$, each having
countably many branches and each branch is a conformal bijection:

\begin{enumerate}\itemsep2mm
\item $\psi_k$: This branch maps $\mathbb C\setminus (-\infty,z_k\rbrack$ to $\Omega_{k}$ for each $k\in \mathbb Z$. Each $\psi_k$ forms a conformal bijection.
  \[
  \psi_{k}:\mathbb C\setminus (-\infty,z_k\rbrack\to \Omega_{k}, \quad
  \text{for } k\in \mathbb Z.
  \]

\item $\widetilde{\psi}_k$: This branch maps $\mathcal A_k$ to $D_{k}$
  for each $k\in \mathbb Z$. Each $\widetilde \psi_k$ forms a
  conformal bijection.
  \[
  \widetilde\psi_{k}:\mathcal A_k\to D_{k}, \quad \text{for } k\in \mathbb Z.
  \]

\item $\widehat{\psi}_k^+$: This branch maps $\mathbb C_{k+}$ to
  $D_{k}$ for each $k\in \mathbb Z$. Each $\widehat\psi_{k}^+$ forms a
  conformal bijection.
  \[
  \widehat\psi_{k}^+:\mathbb C_{k+}\to D_{k}, \quad \text{for } k\in \mathbb Z.
\]

\item $\widehat{\psi}_k^-$: This branch maps $\mathbb C_{k-}$ to
  $D_{k}$ for each $k\in \mathbb Z$. Each $\widehat\psi_{k}^-$ forms a
  conformal bijection.
  \[
  \widehat\psi_{k}^-:\mathbb C_{k-}\to D_{k}, \quad \text{for } k\in \mathbb Z.
  \]
\end{enumerate}
To divide the entire $w$-plane, we add the ``upper'' boundary to the
domain below, analogous to the methodology applied in the cases
$a=\frac{1}{2}$ and $a=\frac{1}{4}$. Consistent with previous cases,
the branch point structure reveals that the branch points are given by
\[
\textbf{BP}_{a}=\textbf{CP}_{a}\cup\{0\}.
\]
The process of constructing the Riemann surface connected to $\psi$
and describing the monodromy group is analogous to the cases
$a=\frac{1}{2}$ and $a=\frac{1}{4}$, and thus, will be omitted here.


\section{Special Cases}\label{sec:specialcases}

We end this article by investigating some special cases of particular
interest. Specifically, we explore the behavior of complex branches of
the function $\psi$ as the variable $a$ approaches $0$, $1$, and
$\frac{1}{3}$.\\

\noindent\emph{Case $a\to 0^+$:}\\

Consider the situation where $a$ is infinitesimally close to zero. We
propose that, under this condition, $\psi$ approximates the classical
Lambert $W$ function, $W$. This can be deduced in two ways. To begin
with, when $a$ approaches $0$, the following approximation holds:
\begin{equation*}
  e^z\sinh(az)=aze^z\frac{\sinh(az)}{az}\approx aze^z,
\end{equation*}
which consequently implies:
\begin{equation*}
  \psi(s)\approx W(s/a).
\end{equation*}

This can be further verified through a careful analysis of the
branches. Furthermore, $x_a\approx \frac{1}{e}$ and for a fixed
neighborhood of zero, no other critical points exist aside from
$z_0\approx -\frac{a}{e}$. This gives us two branch points,
$\textbf{BP}_a=\{-\frac{a}{e},0\}$. In addition, there is a single
branch $\psi_0$ and infinitely many branches $\widetilde\psi_k$.\\

\noindent\emph{Case $a\to 1^-$}:\\

As $a$ converges towards $1$, the complex branches of $\psi$ conform
to the complex branches defined by the logarithmic function. If $a\to
1$ then the function $f(z)$ approaches $e^z\sinh(z)$, implying:
\begin{equation*}
  \psi(s)=\frac{1}{2}\log(2s+1).
\end{equation*}
Considering the Jacobian (\ref{jacobian}), it is always positive and
the imaginary part of $f(z)$ equals zero if and only if $\eta=k\pi$,
where $k\in \mathbb Z$. Additionally, $z_k\to -\frac{1}{2}$, which
will be the sole branch point. Hence, our branches will eventually
become the branches of the logarithmic function.\\

\noindent\emph{Case $a=\frac{1}{3}$}:\\

This scenario stands out as a special case where $\frac{a+1}{a-1}$
belongs to the set of natural numbers, $\mathbb{N}$. For
$a=\frac{1}{3}$, we are able to provide explicit formulas for the
branches as follows:
\begin{align*}
  \psi(z) & =\frac{3}{2}\log \left(\frac{1}{2}+\sqrt{2z+\frac{1}{4}}\right), \\
  \widetilde\psi(z) & =\frac{3}{2}\log \left(\frac{1}{2}-\sqrt{2z+\frac{1}{4}}\right).
\end{align*}
In this particular case, there exists only a single critical point,
denoted as
$\textbf{CP}_{1/3}=\{x_{1/3}\}=\{-\frac{1}{8}\}$, and
two branch points,
$\textbf{BP}_{1/3}=\{-\frac{1}{8},0\}$. Furthermore,
$\mathcal{A}_k=\mathbb{C}\setminus(-\infty,0)$ defines the domain of
$\widetilde{\psi}$, while $\mathbb{C}\setminus(-\infty,-\frac{1}{8})$
defines the domain of $\psi$.


\section{Conclusions and Future Work}\label{sec:case(3)}

Motivated by the intricate structures of the classical Lambert $W$ function's complex branches (see e.g.~\cite{MezoBook}) and its generalizations as discussed by Mez\H{o}~\cite{Mezo21}, this study introduces and examines complex branches of the inverse function $\psi$ of $f(w)=\sinh(aw)e^w$, where the parameter $a$ lies within the range $0<a<1$. The connection of function $\psi$ with $p,q$-binomial coefficients, and the Lenz-Ising model is detailed in \cite{AhagCzyzLundow}.

Our analysis of the complex branches of $\psi$ and their associated Riemann surfaces hinges on the parameter $a$. We categorize our analysis based on whether $\frac{1+a}{1-a}$ belongs to $\mathbb{N}$ or not. The most challenging scenario arises when $a \notin \mathbb{Q}$, a topic that remains unresolved, as noted in the closing remarks of this section.

We explore special limit cases as $a$ approaches $0^+$ and $1^-$. In the former scenario, there is a discernible link between $\psi$ and the classical Lambert $W$ function. We provide an explicit expression for $\psi$ in the case where $a=\frac{1}{3}$.

The analysis becomes significantly more involved for values of $a$ outside the rational numbers, $\mathbb{Q}$. In these cases, the function $f$ exhibits a notable lack of periodicity, making each irrational value of $a$ require a separate examination.

We conjecture that the analysis of complex branches for such $a$ values mirrors the approach employed for rational values, barring the periodicity aspect. Even with the complexity, this analysis should somewhat parallel the earlier scenario, with only subtle variations in how the $w$-plane is partitioned. This plane is segmented by the curves $\Gamma_k$, as initially introduced for $a=\frac{1}{4}$. The positioning of these ``parabolic'' and ``cubic'' curves, denoted as $\Gamma_k$, relies on the parameter $a$. In such scenarios, we suspect to encounter branches of types $\psi_k$, $\widetilde{\psi}_k$, $\widehat{\psi}_k^+$, and $\widehat{\psi}_k^-$.

\section*{Disclosure Statement}

No potential conflict of interest was reported by the authors.


\begin{thebibliography}{99}

\bibitem{AhagCzyzLundow} {\AA}hag P., Czy{\.z} R., Lundow P. H., On a
  generalised Lambert $W$ branch transition function arising from
  $p,q$-binomial coefficients. Appl. Math. Comput. 462 (2024), Paper No. 128347, 20 pp.

\bibitem{BariczMezo} Baricz \'{A}, Mez\H{o} I., On the generalization
  of the Lambert $W$ function. Trans. Amer. Math. Soc. 369 (2017),
  no. 11, 7917-7934.

\bibitem{Beardon} Beardon A. F., The principal branch of the Lambert
  $W$ function. Comput. Methods Funct. Theory 21 (2021), no. 2,
  307-316.

\bibitem{Beardon2} Beardon A. F., Winding numbers, unwinding numbers,
  and the Lambert $W$ function. Comput. Methods Funct. Theory 22
  (2022), no. 1, 115-122.

\bibitem {BhamidiSteeleZaman} Bhamidi S., Steele J. M., Zaman T.,
  Twitter event networks and the superstar model,
  Ann. Appl. Probab. 25(5) (2015), 2462-2502.

\bibitem{CorlessGonnetHareJeffreyKnuth} Corless R. M., Gonnet G. H.,
  Hare D. E. G., Jeffrey D. J., Knuth D. E., On the Lambert $W$
  function. Adv. Comput. Math. 5 (1996), no. 4, 329-359.


\bibitem{Kozlov} Kozlov M., Tulendinova A., Kim J., Ellis
  G., Skrzypacz P., Oscillations of retaining wall subject to Grob’s
  swelling pressure. Scientific Reports 12(1) (2022), p.12224.

\bibitem{LundowRosengren} Lundow P. H., Rosengren A. On the $p,
  q$-binomial distribution and the Ising model. Philos. Mag. 90
  (2010), no. 24, 3313-3353.

\bibitem{Mezo21} Mez\H{o} I., The Riemann surface of the $r$-Lambert
  function. Acta Math. Hungar. 164 (2021), no. 2, 439-450.

\bibitem{MezoBook} Mez\H{o} I., The Lambert $W$ function -- its
  generalizations and applications. Discrete Mathematics and its
  Applications (Boca Raton). CRC Press, Boca Raton, FL, 2022. xxi+252
  pp.

\bibitem{ScottFreconGrotendorst} Scott T. C., Frecon M. A.,
  Grotendorst J., New approach for the electronic energies of the
  hydrogen molecular ion, Chem. Phys., 324 (2006), 323-338.

\bibitem{ScottMann} Scott T. C., Mann R. B., General relativity and
  quantum mechanics: Towards a generalization of the Lambert $W$
  function, Appl. Algebra Engrg.  Comm. Comput. 17 (2006), 41–47.

\bibitem{TrefethenWeideman} Trefethen L. N., Weideman J. A. C., The
  exponentially convergent trapezoidal rule. SIAM Rev. 56 (2014),
  no. 3, 385-458.

\end{thebibliography}
\end{document}